\documentclass[11pt]{article}

\catcode`\@=11
\def\marginnote#1{}

\newcount\hour
\newcount\minute
\newtoks\amorpm
\hour=\time\divide\hour by60 \minute=\time{\multiply\hour by60
\global\advance\minute by-\hour}
\edef\standardtime{{\ifnum\hour<12
\global\amorpm={am}%
        \else\global\amorpm={pm}\advance\hour by-12 \fi
        \ifnum\hour=0 \hour=12 \fi
        \number\hour:\ifnum\minute<10
        0\fi\number\minute\the\amorpm}}
\edef\militarytime{\number\hour:\ifnum\minute<10
0\fi\number\minute}

\def\draftlabel#1{{\@bsphack\if@filesw {\let\thepage\relax
   \xdef\@gtempa{\write\@auxout{\string
      \newlabel{#1}{{\@currentlabel}{\thepage}}}}}\@gtempa
   \if@nobreak \ifvmode\nobreak\fi\fi\fi\@esphack}
        \gdef\@eqnlabel{#1}}
\def\@eqnlabel{}
\def\@vacuum{}
\def\draftmarginnote#1{\marginpar{\raggedright\scriptsize\tt#1}}

\def\draft{\oddsidemargin -.5truein
        \def\@oddfoot{\sl preliminary draft \hfil
        \rm\thepage\hfil\sl\today\quad\militarytime}
        \let\@evenfoot\@oddfoot \overfullrule 3pt
        \let\label=\draftlabel
        \let\marginnote=\draftmarginnote
   \def\@eqnnum{(\theequation)
   \rlap{\kern\marginparsep\tt\@eqnlabel}%
\global\let\@eqnlabel\@vacuum} } %%%%%%%%%%%
\relax

\input{amssym.def}
\input{amssym}
\usepackage{eufrak}

\newtheorem{theorem}{Theorem}

\newtheorem{corollary}[theorem]{Corollary}
\newtheorem{definition}[theorem]{Definition}

\newtheorem{proposition}[theorem]{Proposition}
\newtheorem{remark}[theorem]{Remark}

\def\beq{\begin{equation}}
\def\be{\begin{equation}}
\def\eeq{\end{equation}}
\def\ee{\end{equation}}
\def\K{{\Bbb K}}
\def\C{{\Bbb C}}
\def\R{{\Bbb R}}
\def\B{{\cal B}}
\def\L{{\cal L}}
\def\M{{\cal M}}
\def\N{{\cal N}}

\def\KK{{\cal K}}
\def\Kh{{\cal K}_\h}
\def\D{{\cal D}}
\def\W{{\cal W}}
\def\qqq{q^{-3}}
\def\q{q^{-1}}

\def\qq{q^{-2}}
\def\pa{\partial}
\def\dd{\tilde{\partial}_t}
\def\hhh{\h^{-1}}

\def\uq{U_q(sl(2))}
\def\Uq{U_q(sl(m))}

\def\ot{\otimes}
\def\h{{\hbar}}

\def\vv{V^{\ot 2}}

\def\Mat{{\rm Mat}}

\def\span{{\rm span}}
\def\Sym{{\rm Sym}}

\def\Tr{{\rm Tr}}

\def\De{{\Delta}}
\def\de{{\delta}}
\def\Gr{{\rm Gr}}
\def\Cas{{\rm Cas}}

\def\al{\alpha}
\def\Al{\overline{\alpha}}

\def\hh{{\frak h}}

\def\la{{\lambda}}

\def\L{{\cal L}}

\def\hh{\displaystyle\frac{\h}{2}}
\begin{document}

\makeatletter
\renewcommand{\theequation}{{\thesection}.{\arabic{equation}}}
\@addtoreset{equation}{section}
\makeatother

\title{Braided  Weyl algebras and differential calculus on $U(u(2))$}
\author{
\rule{0pt}{7mm} Dimitri
Gurevich\thanks{gurevich@univ-valenciennes.fr}\\
{\small\it LAMAV, Universit\'e de Valenciennes,
59313 Valenciennes, France}\\
\rule{0pt}{7mm} Pavel Pyatov\thanks{pyatov@theor.jinr.ru}\\
{\small\it
Faculty of Mathematics, NRU HSE, 101000 Moscow, Russia}\\
{\footnotesize\it \&}\\
{\small\it
Bogoliubov Laboratory of Theoretical Physics, JINR, 141980 Dubna, Russia}\\
\rule{0pt}{7mm} Pavel Saponov\thanks{Pavel.Saponov@ihep.ru}\\
{\small\it
Faculty of Mathematics, NRU HSE, 101000 Moscow, Russia}\\
{\footnotesize\it \&}\\
{\small\it Division of Theoretical Physics, IHEP, 142281
Protvino, Russia} }

\maketitle

\begin{abstract}
On any  Reflection  Equation algebra corresponding to a skew-invertible Hecke
symmetry (i.e. a special type solution of the Quantum Yang-Baxter Equation)
we define analogs of the partial derivatives. Together with elements of the initial
Reflection Equation algebra they generate a "braided analog" of the Weyl
algebra. When  $q\to 1$,  the braided Weyl algebra corresponding to the
Quantum Group $U_q(sl(2))$ goes to the Weyl algebra defined on
the algebra $\Sym((u(2))$ or  that $U(u(2))$ depending on the way of passing to the limit.
Thus, we define partial derivatives on the algebra  $U(u(2))$, find their "eigenfunctions",
and introduce an analog of the Laplace operator on this algebra. Also, we
define the "radial part" of this operator, express it in terms of "quantum
eigenvalues", and sketch an analog of the de Rham complex on the algebra $U(u(2))$.
Eventual applications of our approach are discussed.
\end{abstract}

{\bf AMS Mathematics Subject Classification, 2010:} 17B37, 81R60


\section{Introduction}

Since creation of the Quantum Group (QG) theory numerous attempts of
developing a quantum version of differential calculus were undertaken. This study
was initiated in \cite{WZ} where the role of a quantum function space was played
by a "q-symmetric" algebra of the fundamental  space $V$ equipped with an
action of the QG $\Uq $, and  \cite{W} where this role was played by a compact
matrix pseudogroup.

Essentially, such a pseudogroup is the famous RTT algebra (see \cite{FRT})
associated with a given braiding, i.e. an invertible operator $R:\vv\to\vv $,
satisfying  the Quantum Yang-Baxter Equation
 $$
 (R\ot I)\,(I\ot R)\,(R\ot I)=(I\ot R)\,(R\ot I)\,(I\ot R).
 $$
Hereafter $V$ is a vector space over the ground field $\K$ ($=\R$  or  $\C $) and
$I$ is the identity operator.

In \cite{W} a general scheme of defining differential forms and the de Rham
complex on  a matrix pseudogroup
was suggested. Also, the author considered analogs of
vector fields, introducing them by duality. In some subsequent publications (see
f.e. \cite{IP, FP}) the algebra generated by such fields in the case related to the QG $U_q(sl(m))$
was identified  as a (modified)
Reflection Equation (RE) algebra\footnote{In the paper \cite{W} the deformation
property of algebra in question was disregarded. Note, that only for the $A_n$
series it is possible to construct a quantum differential algebra with good
deformation property. Moreover, in the family of classical simple Lie algebras only for $g=sl(m)$
(and consequently, $gl(m)$) there exist $U_q(g)$-covariant deformations of the
algebras $\Sym(g)$ and $\bigwedge(g)$ with classical dimensions of
homogenous components. $U_q(sl(m))$-covariant deformations of the algebra $\Sym(gl(m))$ and $\bigwedge(gl(m))$
can be constructed with the use of some idempotents
playing the role of symmetrizers (resp.,  skew-symmetrizers) (see \cite{GPS2}).
Whereas the operators used  in \cite{W}  (as well as in a number of papers
devoted to the so-called Woronowicz-Nichols algebras) in a similar construction
are not idempotents and are not motivated by the algebraic structure of the initial
braiding.} (see below).

Recently, we have  (partially\footnote{We disregarded quantum analogs of
differential forms.}) generalized  this construction replacing the RTT algebra by
other quantum matrix algebras, in particular, by an RE algebra. Namely, taking a
copy of the RE algebra (denoted $\M $) as a quantum function algebra, we
treated another copy of this algebra (denoted $\L $) as an analog of the
one-sided or adjoint differential operators. Below, we deal with the total algebra
(denoted $\B(\L,\M)$) where we assume that $\M$ is equipped with a left (right-
invariant) action of $\L$.

Besides, the braiding $R$ coming in the definition of the algebra $\B(\L,\M)$ is
taken to be a Hecke symmetry. This means that it  is subject to the equation
$$
(R-q\, I)\, (R+\q\,I)=0,
$$
where  $q\in\K$ is assumed to be generic. In particular,  such a Hecke symmetry
comes from the  QG $\Uq $. This Hecke symmetry and all related objects will be
called {\em standard}. The RTT and RE algebras associated with a standard
Hecke symmetry are deformations of the algebra $\Sym(gl(m))$. With
the use of other Hecke symmetries we can get analogous deformations of the
super-algebra $\Sym(gl(m|n))$.

Note that any RE algebra associated with a Hecke symmetry has another basis in
which the permutation relations between basic elements are quadratic-linear. So, the RE algebra in this basis
becomes more similar to the enveloping algebra
$U(gl(m|n))$. We call this quadratic-linear algebra
the {\em modified} RE (mRE) algebra\footnote{Note that the
passage to this new basis (or in other words, the isomorphism between the
RE algebra and its modified form) fails for $q=1$. Besides, in order to treat any enveloping algebras as a deformation  of the corresponding symmetric one
 we introduce another  parameter $\h$ in the defining relations of the enveloping algebra (and also in the related mRE algebra).}. Namely in this form
  the RE algebra comes in
constructing quantum analogs of vector fields.

Besides, in \cite{GPS2} we have constructed a representation  category of a mRE algebra looking like that for the enveloping algebra
$U(gl(m|n))$. Moreover, according to our construction   the former category turns into the later one, provided that
the mRE algebra goes to $U(gl(m|n))$ in the limit $q\to 1$.

One of the aim of the present paper is to define a braided counterpart of the partial
derivatives on the RE algebras (modified or not). First, following \cite{GPS3}
we equip an RE algebra $\M $  with a left (i.e. right invariant) action of a mRE algebra $\KK$.
(The both algebras are defined via the same Hecke symmetry.)

Second, by combining the generating matrices
of these two algebras we construct a matrix of "partial derivatives" on $\M$. In the standard case
these "derivatives" turn into the usual ones on the algebra  $\Sym(gl(m))$ as $q\to 1$.
The total algebra  generated by the algebra $\M$ and by these partial derivatives on it, is called the
{\em braided Weyl  algebra} and is denoted $\W(\M)$.

Upon replacing the algebra $\M$ by its modified
form $\N$, we get (after a slight rescaling) the partial derivatives on the algebra
 $\N$. The corresponding Weyl algebra is denoted $\W(\N)$.

Observe that in the case related to the QG $\uq$ the algebra $\W(\M)$ appeared in
\cite{OSWZ,K,AKR} in  construction of a q-analog of the Minkowski space. As for
the braided Weyl algebra $\W(\N)$, it is, up to our knowledge, an absolutely  new
object which has a very interesting limit as $q\to 1$. Namely, assuming the initial
Hecke symmetry $R$ to be a deformation of a super-flip, we get in the limit
$q\rightarrow 1$ the partial derivatives on the algebra $U(gl(m|n))$ and consequently the
Weyl algebra $\W(U(gl(m|n)))$.

Naturally, the usual Leibniz rule valid for the derivatives on a commutative algebra is
modified when we pass to one of the mentioned algebras. On the algebra $U(gl(m))$ this Leibniz rule  can be expressed via the following coproduct
\be
\De(\pa_{n_i^j})=\pa_{n_i^j}\ot 1+1\ot \pa_{n_i^j}+\sum_k \pa_{n_i^k}\ot \pa_{n_k^j},
\label{lei}
\ee
where $n_i^j$ are entries of the generating
matrix of this algebra. This form of the Leibniz rule was first found
by S.~Meljanac and Z.~\v{S}koda.  As for the
braided Weyl algebras $\W(\M)$ (resp.,  $\W(\N)$) such a simple formula for
partial derivatives is not yet found, and their action on the algebra $\M$ (resp.,
$\N$) is defined via the permutations relations between elements of this
algebra and derivatives completed via a counit.

We consider the Weyl algebra $\W(U(gl(2)))$ in detail. More precisely, we pass
to the compact form of the Lie algebra $gl(2,\C)$ and deal with the Weyl algebra
$\W(U(u(2)))$. This form is more convenient for defining  wave operators in our
noncommutative (NC) algebra setting.
Namely, they can be introduced in the classical way but with a new meaning of the
derivatives.

Also, we describe a way of defining the radial part $\De_{rad}$ of the Laplace
operator $\De$. The crucial role in our construction is played by eigenvalues of
the generating matrix of the algebra  $U(u(2))$ (see (\ref{matrN})). These
eigenvalues are defined as roots of the Cayley-Hamilton identity for the
generating matrix and are treated as elements of the algebraic extension
of the center of this algebra.

Our final
formula shows that the operator $\De_{rad}$ expressed via  symmetric function
of these eigenvalues is a second order difference operator. We consider this operator
as a first step in the direction of constructing similar analogs of Calogero-Moser
operators and of radial parts of Laplace operators on super-algebra in the spirit of \cite{B}.
A plan of applying our method to this end is exhibited at the end of the paper.

The paper is organized as follows. In the next section we present a construction
of the braided Weyl algebras on the RE algebras (modified or not) and find their
$q\to 1$ limits. In section 3 we consider a two-dimensional example in detail. As
a result we get  the Weyl algebra $\W(U(u(2)))$. Also, in section 3 we calculate
"eigenfunctions" of the partial derivatives and present an analog of the de Rham
complex on the enveloping algebra $U(u(2))$. Besides, in sections 2 and 3 we
calculate Poisson counterparts arising from deformations in question. In section 4
we define the Laplacian on the algebra $U(u(2))$ and compute its radial part.

{\bf Acknowledgement}
This work was partially supported by the joint RFBR and CNRS grant
09-01-93107-NCNIL-a. The work of P.P. and P.S. was partially supported by the
RFBR grant 11-01-00980-a and by the Higher School of Economics Academic
Fund grant 11-09-0038. Also, the work of P.P. and P.S. was partially supported by the
Higher School of Economics Academic Fund grants 10-01-0013 and  11-01-0042
respectively. D.G. is thankful to S.~Meljanac and Z.~\v{S}koda for valuable
discussions and  the franco-croatian cooperation programm  Egide PHC Cogito 24829NH
for a financial support of his visit to Zagreb  University.

\section{Braided Weyl algebras and their $q\to 1$ limits}

In what follows $R:\vv\to\vv $ is assumed to be a skew-invertible Hecke symmetry
(see \cite{IOP, GPS2} for definitions). The space $V$ is called {\em basic}.

With any skew-invertible Hecke symmetry $R:\vv\to \vv $ we associate a quantum
matrix algebra which is referred to as the Reflection Equation (RE) algebra.
\begin{definition}
The RE algebra is a unital associative algebra over the field $\K$ generated
by entries $l_i^j$ of a matrix $L=\|l_i^j\|$, $1\leq i,j \leq \dim \, V$, which are
subject to the system of commutation relations
\be
R\, L_1 R\, L_1=L_1 R\, L_1 R,\qquad L_1=L \ot I.
\label{RE}
\ee
\end{definition}
In this definition $I$ stands for the unit matrix. Also, we assume  a
basis $\{x_i\}$ to be fixed in the space $V$. Then, in the basis $\{x_i\ot x_j\}$ of
$\vv $ the braiding $R$ is represented by a matrix $\|R_{ij}^{kl}\|$:
$$
R(x_i\ot x_j)=R_{ij}^{kl} x_k\ot x_l.
$$
Hereafter a summation over the repeated indices is assumed. The RE algebra
defined above is a particular cases of more general construction of a quantum
matrix algebra introduced in \cite{IOP}.

Let us consider two copies of the RE algebra: one of them, generated by entries
of $M=\|m^i_j\|$ and denoted  $\M $, plays the role of the quantized function algebra
$\Sym(gl(m))$. The other one, generated by entries of $L=\|l^i_j\|$ and denoted  $\L $,
plays the role of quantized right-invariant differential operators on $\M$. Their action
is encoded in the following permutation relations between two families of generators
\be
R\, L_1\, R\, M_1=M_1\,R\, L_1\,R^{-1}. \label{perm}
\ee
Thus, the whole algebra $\B(\M,\L)$ is generated by entries of two matrices $M$
and $L$ which satisfy the RE algebra permutation relations (\ref{RE}) and, additionally,
subject to (\ref{perm}).

Permutation relation (\ref{perm}) defines one of the possible doubles of the two
reflection equation algebras. This double was intensively studied as an algebra
possessing coproduct and coaddition structures (see \cite{IV}
and references therein) and also in application to a q-deformation of the Poincare
algebra \cite{K,AKR}.

To get the action of an element $l\in \L$ to an element $m\in \M$ we proceed
as follows. Note, that the permutation relations enable us to reduce any element
of the whole algebra $\B(\M,\L)$ to that of the tensor product $\M\ot \L$. So, given
an element $l\ot m\in \L\ot \M$, we reduce it to the form of an element
from $\M\ot \L$, and then apply a counit $\varepsilon: \L\rightarrow \K$ to the
components from $\L$.  Here $\varepsilon$ is an ingredient of a braided bi-algebra
structure on any RE algebra which was discovered by Sh.~Majid \cite{M}.
On generators of RE algebra the coproduct and counit of the braided bi-algebra are
defined as follows
$$
\De(1_\L)=1_\L \ot 1_\L,\quad \De(l_i^j)=l_i^k \ot l_k^j,\quad \varepsilon(1_\L)=
1,\quad \varepsilon(l_i^j)=\de_i^j.
$$
(Hereafter, $1_{\cal A}$ stands for the unit of a unital algebra ${\cal A}$.)

Thus, in order to compute the action $l_i^j(m)$ we employ the  following chain of
transformations
\be
l_i^j(m) \equiv l_i^j( m\, 1_{\M})\stackrel{(\ref{perm})}
{\longrightarrow}m'\, l_{i'}^{j'}(1_\M)\stackrel{\mbox{\tiny def}}{=}
m'\,\varepsilon(l_{i'}^{j'})=m'\, \de_{i'}^{j'},
\label{action-ch}
\ee
where we set by definition $l(1_\M)=\varepsilon(l)1_\M$.

This action can be extended to all polynomials in generators $l_i^j$ in a straightforward
way and putting additionally $1_\L(m)=m$ for any $m\in \M $, we  define an action
of any polynomial in $l_i^j$ onto the whole algebra $\M$.

\begin{remark} \label{rem3} \rm
However, we can renormalize the action of the generators $l_i^j$ on $1_{\M}$ by
putting $l_i^j(1_{\M})=\la\,\de_{i}^{j}\, 1_{\M} $ where $\la$ is a nontrivial numerical factor.
Thus, we get a representation of the algebra $\L$ multiple to the previous one.
Observe that the algebra $\L$ allows an automorphism $l_i^j\to \la\,  l_i^j$.
\end{remark}

Now, let us define  a modified version of the RE algebra.

\begin{definition}
Let $R: \vv\to\vv $ be a Hecke symmetry. The modified RE algebra $\Kh $ is a
unital associative algebra generated by entries $k_i^j$ of the matrix $K=\|k_i^j\|$,
$1\le i,j\le \dim \, V$, which are subject to the system of commutation relations
\be
R\, K_1\, R\, K_1 - K_1\,R\, K_1\,R=\h\,(R\, K_1-K_1\, R),\qquad K_1=K\otimes
I,\quad \h\in \K.
\label{mREA}
\ee
\end{definition}
For the case $\h=1$ we shall omit the subscript in notation $\Kh$.

Unless $q=\pm 1$, the algebras $\L$ and $\Kh$ are isomorphic. The
isomorphism is explicitly defined by the following relation between two sets of
generators\footnote{Since a rescaling  $L\to \lambda\, L,$ $\lambda\in \K,$
$\lambda\not=0$ is an automorphism of the algebra $\L$ the map (\ref{iso})
can be also rescaled. Thus, in \cite{GS2} we used another normalization of this map.}
\be
1_\L= 1_\KK,\quad L= \h\,1_{\KK} I\, -(q-\q)\,K.
\label{iso}
\ee
Nevertheless, the limits of the algebras $\L$ and $\Kh$ as $q\to 1$ are not
isomorphic to each other since the map (\ref{iso}) degenerates.

Note that the above coproduct in the RE algebra being rewritten in the generators
$k_i^j$  takes the form (hereafter we set $\h=1$)
\be
\De(1_\KK)=1_\KK\ot 1_\KK,\quad \De(k_i^j) = k_i^j\ot 1_\KK+1_\KK\ot
k_i^j-(q-\q) \,k_i^s\ot k_s^j.
\label{cco}
\ee
A way of extending this coproduct to the whole algebra $\KK $ is described in
\cite{GPS2}.

The last term  of this formula disappears at $q\to 1$ and the coproduct turns into
the usual one coming in the Hopf structure of the algebra $U(gl(m|n))$.
As for the counit we choose it in the form
\be
\varepsilon(1_\KK)=1,\quad  \varepsilon(k_i^j)=0.
\label{K-counit}
\ee

So, we rewrite the defining relations of the algebra $\B(\M,\L)$ in terms of  the
new generators $k_i^j$ (see (\ref{iso})) and denote this algebra as $\B(\M,\KK)$.
The algebra $\B(\M,\KK)$ contains an RE subalgebra $\M $, a modified RE
subalgebra $\KK $, while the permutation relation between the generating
matrices of these subalgebras  becomes
\be
R\, K_1 R\, M_1=M_1 R\, K_1  R^{-1}+ R\, M_1.
\label{km-rel}
\ee

By the same method as above
we can define an action of  the algebra $\KK$ on $\M$: for this end we employ the
permutation relations (\ref{km-rel}) and the counit. Thus, we get the braided counterparts
of usual right-invariant vector fields. Namely, as $q\to 1$ the elements $k_i^j$ treated
as operators turn into such vector fields (see \cite{GPS3} for more detail).

Now, consider a matrix $D=M^{-1}\, K$. The matrix  $M^{-1}$ inverse to $M$ can be
defined  via the Cayley-Hamilton (CH) identity established in \cite{GPS1} for the
generating matrix $M$ of any RE algebra $\M $. For the existence of $M^{-1}$ we
should only assume the lowest coefficient (which is always central) of the CH identity
to be invertible. More precisely, this operation can be realized as a localization by the
mentioned central element.

\begin{remark} \rm Emphasize that in the frameworks of the method suggested in
\cite{GR} of inverting matrices with entries from a NC algebra one requires the
invertibility of a large number of elements of this algebra. In contrast with \cite{GR},
our method can be applied only to the generating matrix of the RE algebra. Say,
other matrices with entries from this algebra are not subject to a similar CH identity
and our method of inversion is not valid for them.
\end{remark}

Following the classical pattern we treat entries of the matrix $ D=\|\partial_i^j\|$ as
analogs of partial derivatives (momenta) : $\partial_i^j=\partial/\partial{m_j^i}$.
Moreover, we can give a literal sense to this treatment, by defining an action
of the partial derivatives onto the algebra $\M $. To this end we  explicitly present
the permutation relations between the  generators $m_i^j$ and $\partial_k^l$.

\begin{proposition}
\label{prop:4}
{\rm
Consider a unital associative algebra $\W(\M)$ generated by entries of the matrices
$M$ and $D$. Then the defining relations of the algebra ${\cal B}(\M,\KK)$ leads to the
following relations between generators of $\W(\M)$
\be
\begin{array}{rcl}
R\,M_1R\,M_1&=&M_1R\,M_1R,\\
R^{-1}D_1R^{-1}D_1&=&D_1R^{-1}D_1R^{-1},\\
D_1 R\,M_1 R&=&R\,M_1R^{-1}D_1+R.
\end{array}
\label{weyl}
\ee
}
\end{proposition}

{\bf Proof.} We only should check the second an third relations. For this purpose
we need some additional formulae for commutation relations of entries of
$M^{-1}$ with each other and with entries of $K$. First, the RE for $M$ directly
leads to the relation
$$
M^{-1}_1R^{-1}M_1^{-1}R^{-1} = R^{-1}M_1^{-1}R^{-1}M_1^{-1}.
$$
Then, we rewrite the defining relation (\ref{km-rel}) in an equivalent form
$$
K_1R^{-1}M_1^{-1} = R^{-1}M_1^{-1}RK_1R - R^{-1}M_1^{-1}R.
$$
Now, we make the identical transformations:
\begin{eqnarray*}
R^{-1}D_1R^{-1}D_1 &=& R^{-1}M_1^{-1}\underline{K_1R^{-1}M_1^{-1}}K_1 =
\underline{R^{-1}M_1^{-1}R^{-1}M_1^{-1}}RK_1RK_1\\
& - &\underline{R^{-1}M_1^{-1}R^{-1}M_1^{-1}}RK_1
=M_1^{-1}R^{-1}M_1^{-1}K_1RK_1 - M_1^{-1}R^{-1}M_1^{-1}K_1.
\end{eqnarray*}
Here we underline the part of an expression which undergoes an identical
transformation at the next step of calculations.

On the other hand we get
\begin{eqnarray*}
D_1R^{-1}D_1R^{-1} &=& M_1^{-1}\underline{K_1R^{-1}M_1^{-1}}K_1R^{-1} =
M_1^{-1}R^{-1}M_1^{-1}\underline{RK_1RK_1}R^{-1} - \\
M_1^{-1}R^{-1}M_1^{-1}
RK_1R^{-1}
&=& M_1^{-1}R^{-1}M_1^{-1}K_1RK_1 - M_1^{-1}R^{-1}M_1^{-1}K_1.
\end{eqnarray*}
Here we used the defining relations of the modified RE algebra $\KK $ in order to
transform the term $RK_1RK_1$ in the last step of transformations.  Since the
right hand sides of the above equalities coincide, we conclude that permutation rules
for $D$ is indeed as is claimed in (\ref{weyl}).

To prove the last line relation of the system (\ref{weyl}) we use the following
auxiliary formulae
$$
M^{-1}_1R^{-1}M_1R = RM_1R^{-1}M_1^{-1}\quad {\rm and}\quad
K_1RM_1 = R^{-1}M_1RK_1R^{-1} + M_1.
$$
Now, we have
\begin{eqnarray*}
D_1RM_1R &=& M_1^{-1}\underline{K_1RM_1}R = \underline{M_1^{-1}R^{-1}
M_1R}K_1 + R\\
&=&RM_1R^{-1}M_1^{-1}K_1 +R = RM_1R^{-1}D_1 +R.
\end{eqnarray*}
The proof is completed.\hfill\rule{6pt}{6pt}

\begin{definition} We call the algebra ${\cal W}(\M)$ defined in (\ref{weyl}) a
braided Weyl algebra.
\end{definition}
Note, that the subalgebra $\D\subset \W(\M)$ generated by $\partial_i^j$ is also an RE
algebra associated with the Hecke symmetry $R^{-1}$.

As we noticed in Introduction, in the case related to the QG $\uq$ the relations
(\ref{weyl}) appeared in \cite{OSWZ} and \cite{K,AKR} in a study of the q-Minkowski space.

Now, we define an action of partial derivatives onto elements of the algebra $\M$
via the same method as above with the help of an additional requirement
\be
\partial_i^j(1_{\cal W}) = 0
\label{cu}
\ee
which is a direct consequence of the relation $D=M^{-1}K$ and counit (\ref{K-counit}).

Our next aim is to define similar derivatives on the mRE algebra. Above we noticed that
any RE algebra admits a rescaling automorphism. As a consequence, we can get any
nonzero factor at the summand $R$ in the last line of the system (\ref{weyl}). Let us choose
this factor to be equal $-(q-\q)$ and pass from the RE algebra $\M $ to its modified version
with generating matrix $N=\|n_i^j\|$ by applying a shift similar to (\ref{iso})
\be
1_\N=1_\M,\,\,\,M=\h\,1_\M I-(q-\q)\, N.
\label{shift}
\ee
In this way we get a braided Weyl $\W(\N_\h)$  (bellow we omit the subscript $\h$)
defined by the following
relations
\be
\begin{array}{rcc}
R\,N_1R\,N_1-N_1R\,N_1R&=&\h\,(R\,N_1-N_1R),\\
R^{-1}D_1R^{-1}D_1&=&D_1R^{-1}D_1R^{-1},\\
D_1R\,N_1R-R\,N_1R^{-1}D_1&=&R+\h\,D_1R.
\end{array}
\label{set}
\ee

According to our scheme, in order to treat entries of the matrix $D$ as operators
we should complete the above permutation relations with an action of partial derivatives
on the unit element of the algebra. Let us define it by the same formula (\ref{cu}).
The action of the derivatives onto generators $n_k^l$ is the same as on $m_k^l$.

Now, assume a Hecke symmetry $R$ to be a deformation of a super-flip $P$:
$R\stackrel{q\rightarrow 1} {\longrightarrow } P$. Here $P$ acts on the
tensor square of a superspace $V$:
$$
P : \vv\to \vv,\qquad  V=V_0\oplus V_1,\quad \dim V_0=m,\quad \dim
V_1=n,
$$
where $V_0$ (resp., $V_1$) is the even (resp., odd) component of the space
$V$. Then, passing in (\ref{set}) to the limit $q\to 1$ we get a Weyl algebra $\W(U(gl(m|n)_\h))$
on a NC algebra $U(gl(m|n)_\h)$\footnote{Though in \cite{GPS3} we dealt with even
Hecke symmetries, all results of that paper are valid in the general case.}. The
defining relations of this algebra are
\be
\begin{array}{rcc}
P\,N_1P\,N_1-N_1P\,N_1P&=&\h\,(P\,N_1- N_1P), \\
P\,D_1P\,D_1&=&D_1P\,D_1P,\\
D_1P\,N_1 P-P\,N_1P\,D_1&=&P+\h\,D_1P.
\end{array}
\label{sett}
\ee
The first line of the above system  is just the defining relations of the algebra
$U(gl(m|n)_\h)$ (but with a basis slightly different from the usual one formed by the matrix units).
The second one means that the partial derivatives form a
super-commutative algebra and  the third line exhibits the permutation relations
between two  ingredients of the algebra  $\W(U(gl(m|n)_\h))$. Note that the
algebra $\W(U(gl(m|n)_\h))$ is a one-parameter deformation of the super-Weyl
algebra $\W(\Sym(gl(m|n)))$, whereas the algebra $\W(\N)$ above is a
two-parameter deformation of the same algebra (the parameters are $\h$ and $q$).

Now, discuss the Poisson  counterparts of the algebra $\W(U(gl(m)))$ (i.e. we
restrict ourselves to the even case $n=0$). Let $\W(\M)_0$ be the algebra defined by (\ref{weyl})
but without the last term $R$ in the third relation. This algebra is a graded quadratic one,
i.e. it is defined by quadratic relations on generators. For a generic $q$ the dimensions
of its homogeneous components are classical (i.e. they equal those of the space
$\Sym(gl(m))^{\ot 2}$). It is so since this algebra equals (as a set) the tensor product of
two RE algebras and for them this property was proved in \cite{GPS2}.

Besides, consider the quadratic-linear  algebra $\W(\N)_0$ which is also defined
by canceling the term $R$ in the third defining relations (\ref{set}) of the algebra
$\W(\N)$. Since the passage from the algebra $\W(\M)_0$ to the algebra $\W(\N)_0$
can be done via (\ref{shift}), it is easy to see that the graded algebra $\Gr\,\W(\N)_0$
associated with $\W(\N)_0$ is isomorphic to $\W(\M)_0$. This entails the
existence of a Poisson pencil on the commutative algebra which is the $q\to 1$
limit of the algebra $\W(\M)_0$. In the next section we explicitly write down these
Poisson brackets in a two-dimensional example.

\begin{remark}\rm
Observe  that if $n=0$ and  $\h=0$ the  permutation relations from
(\ref{sett})
completed with the  counit are equivalent to the usual Leibniz rule.
It is interesting to find a similar rule in a general case.
In the case $n=0$ and  $\h\not=0$ the Leibniz rule takes the form
\be
\De(\pa_i^j)=\pa_i^j\ot 1+1\ot \pa_i^j+\h\, \pa_k^j\ot \pa_i^k
\label{Lf}
\ee
found S.~Meljanac and Z.~\v{S}koda. (For $\h=1$ this formula turns into (\ref{lei}) since $\pa_i^j=\pa_{n^i_j}$.)

This form together with the rules $\pa_i^j(1_\W)=0$ and $\pa_i^j
(n_k^l)=\de_k^j\, \de_i^l$ enables us to compute the action of partial derivatives
on any element from
$U(gl(m)_\h)$ via the usual formula
$$\pa_i^j(f\cdot g)=\cdot \De(\pa_i^j)(f\ot g)$$
where $\cdot$ stand for the product in the algebra $U(gl(m)_\h)$.

Emphasize the similarity of the coproducts   (\ref{cco}) and (\ref{Lf}).
\end{remark}

It is also interesting to compare the Leibniz rule (\ref{Lf}) and that valid for the
momenta on the $\kappa$-Minkowski space (see \cite{MR}).

\section{Example: $m=2, n=0$}

Let $V$ be a two-dimensional vector space with a fixed basis $\{x,y\}$. Consider
a Hecke
symmetry, which in the basis  $\{x\ot x,\, x\ot y,\, y\ot x,\, y\ot y\}$ of $\vv $ is
represented
by the matrix
$$ \left(\begin{array}{cccc}
q&0&0&0\\
0&q-\q&1&0\\
0&1&0&0\\
0&0&0&q
\end{array}\right).$$
This matrix is just the  product of the usual flip and the image  of the universal
$\uq $ R-matrix in the space $\vv $. Let
\be
N=\left(\begin{array}{cc}
n_1^1&n_1^2\\
n_2^1&n_2^2
\end{array}\right)=\left(\begin{array}{cc}
a&b\\
c&d
\end{array}
\right)
\label{MMM}
\ee
be the generating matrix of the modified RE algebra $\N $.
The  commutation relations (\ref{mREA}) in this case
are as follows
$$
\begin{array}{ll}
q ab-\qq ba=\h b\qquad & q(bc-cb)=((q-q^{-1})\, a-\h)(d-a)\\
\rule{0pt}{4.5mm}
q ca-\qq ac=\h c &  q(db-bd)=((q-q^{-1})\,a -\h)b\\
\rule{0pt}{4.5mm}
ad-da=0 &   q(cd-dc)=c((q-q^{-1})\, a -\h)\\
\end{array}
$$
where we omit the symbol of the unit element  in the linear
combination $(q-q^{-1})\,a-\h $.

Upon setting $\h=0$ we get the defining relations of the corresponding non-
modified
RE
algebra $\M $. We keep the same letters for  entries of the corresponding
generating matrix $M$.

Turn now to the algebra $\D $, which is generated by  entries of the matrix
\be
D=\left(\begin{array}{cc}
\pa_a&\pa_c\\
\pa_b&\pa_d
\end{array}\right).
\label{DDD}
\ee
The relations between these partial derivatives  are the same in both algebras
$\W(\M)$ and
$\W(\N)$.
Their explicit form is
\be
\begin{array}{ll}
\pa_a\,\pa_b-\pa_b\,\pa_a=-(q^2-1) \,\pa_d\,\pa_b\qquad &
\pa_b\,\pa_c-\pa_c\,\pa_b=(q^2-1) \,(\pa_d-\pa_a)\,\pa_d\\
\rule{0pt}{4.5mm}
\pa_a\,\pa_c-\pa_c\,\pa_a=(q^2-1) \,\pa_c\,\pa_d &
q^2\,\pa_b\,\pa_d-\pa_d\, \pa_b=0\\
\rule{0pt}{4.5mm}
\pa_a\, \pa_d-\pa_d\, \pa_a= 0 &
\pa_c\, \pa_d-q^2\,\pa_d\, \pa_c=0.
\end{array}
\label{DDDD}
\ee

The permutation relations between derivatives and generators of the algebra $\N
$
are as follows
$$
\begin{array}{l}
\pa_a\, a =q^{-1}+q^{-2}a\, \pa_a-(1-q^{-2})\,b\, \pa_b +q^{-1}\h\,\pa_a \\
\rule{0pt}{4.5mm}
\pa_a \,b =b\, \pa_a-(1-q^{-2})\,a\, \pa_c+(q-q^{-1})^2\,b\,\pa_d +q^{-1}\h\,\pa_c\\
\rule{0pt}{4.5mm}
\pa_a \,c =q^{-2}\, c\, \pa_a+(1-q^{-2})\,(a-d)\, \pa_b\\
\rule{0pt}{4.5mm}
\pa_a \,d =d\, \pa_a+(1-q^{-2})\,(b\, \pa_b-c\,\pa_c) -(q-\q)^2\, (a-d)\,\pa_d+q(1-
q^{-2})^2
+q(1-q^{-2})^2\h\,\pa_a\\
\rule{0pt}{6mm}
\pa_b\, a=q^{-2}a\,\pa_b+q^{-1}\h\, \pa_b  \\
\rule{0pt}{4.5mm}
\pa_b \,b=q^{-1}+q^{-2} b\, \pa_b-(1-q^{-2})\, a\,\pa_d+q^{-1}\h\, \pa_d \\
\rule{0pt}{4.5mm}
\pa_b\, c=c\, \pa_b \\
\rule{0pt}{4.5mm}
\pa_b\, d=d\,\pa_b-(q^2-1)\, c\, \pa_d\\
\rule{0pt}{6mm}
\pa_c \,a=a\, \pa_c-(q^2-1)\, b\,\pa_d \\
\rule{0pt}{4.5mm}
\pa_c \,b=b\, \pa_c\\
\rule{0pt}{4.5mm}
\pa_c \,c=q^{-3}+q^{-2}c\, \pa_c+(1-q^{-2})(a-d)\,\pa_d-(q^{-2}-q^{-4})\,a\,\pa_a
+(1-q^{-2})^2\,b\,\pa_b+ q^{-3}\h\,\pa_a\\
\rule{0pt}{4.5mm}
\pa_c \,d=q^{-2}d\, \pa_c+(1-q^{-2})\,(2-q^2)\,b\,\pa_d-(1-q^{-2})\,b\,\pa_a
+ (1-q^{-2})^2\,a\, \pa_c + q^{-3}\h\,\pa_c\\
\rule{0pt}{6mm}
\pa_d\, a=a\,\pa_d\\
\rule{0pt}{4.5mm}
\pa_d\, b=q^{-2} b\, \pa_d\\
\rule{0pt}{4.5mm}
\pa_d \,c= c\,\pa_d-(q^{-2}-q^{-4})\, a\, \pa_b+q^{-3}\h\,\pa_b    \\
\rule{0pt}{4.5mm}
\pa_d\, d=q^{-3}+q^{-2}\, d\,\pa_d+(1-q^{-2})^2\,a\, \pa_d-(q^{-2}-q^{-4})\,b\,
\pa_b+q^{-3}\h\,\pa_d
\end{array}
$$

Thus, the result of applying the partial derivatives to the generators of the algebra
$\W(\N)$ is (we only exhibit nontrivial terms)
$$\pa_a(a)=\q,\,\,\pa_b(b)=\q,\,\,\pa_c(c)=\qqq,\,\,\pa_d(d)=\qqq. $$

This system becomes much more simple in the limit $q\to 1$ :
$$
\begin{array}{ll}
\pa_a\,a-a\, \pa_a=1+\h\, \pa_a\qquad&\pa_a\,b-b\, \pa_a=\h\, \pa_c\\
\rule{0pt}{4.5mm}
\pa_a\,c-c\,\pa_a=0&\pa_a\,d-d\, \pa_a=0\\
\rule{0pt}{6mm}
\pa_b\,a-a\, \pa_b=\h\, \pa_b & \pa_b\,b-b\, \pa_b=1+ \h\, \pa_d\\
\rule{0pt}{4.5mm}
\pa_b\,c-c\,\pa_b=0 & \pa_b\,d-d\, \pa_b=0\\
\rule{0pt}{6mm}
\pa_c\,a-a\, \pa_c=0& \pa_c\,b-b\, \pa_c=0\\
\rule{0pt}{4.5mm}
\pa_c\,c-c\, \pa_c=1+\h\,\pa_a &\pa_c\,d-d\, \pa_c=\h\,\pa_c\\
\rule{0pt}{6mm}
\pa_d\,a-a\, \pa_d=0 & \pa_d\,b-b\, \pa_d=0\\
\rule{0pt}{4.5mm}
\pa_d\,c-c\, \pa_d=\h\, \pa_b & \pa_d\,d-d\, \pa_d=1+ \h\,\pa_d.
\end{array}
$$

These are the permutation relations of the partial derivatives with generators of
the algebra $U(gl(2)_\h)$, which are subject to relations
$$
[a,b]=\h\, b,\quad[a,c]=-\h\, c,\quad [a,d]=0,\quad [b,c]=\h\,(a-d),\quad
[b,d]=\h\, b,\quad [c,d]=-\h\,c.
$$
Also, it follows from (\ref{DDDD}) that all partial derivatives commute with each
other. Consequently, they  form a commutative subalgebra $\D $ of the algebra
$\W(U(gl(2)_\h))$.

It is convenient to pass from $gl(2)_\h $ to the compact form $u(2)_\h $.
Introducing new generators as follows
$$
t={{1}\over{2}}(a+d),\quad     x={{i}\over{2}}(b+c),\quad y= {{1}\over{2}}(c-b),\quad
z={{i}
\over{2}}(a-d)
$$
we get the defining relations of the algebra $u(2)_\h $
$$
[x, \, y]=\h\, z,\quad [y, \, z]=\h\, x,\quad[z, \, x]=\h\, y,\quad [t, \, x]=[t, \, y]=[t, \, z]=0.
$$

The change of generators in the partial derivatives is usual:  $\partial_t =
\pa_a+\pa_d$, etc. However, for the future convenience we prefer using the
"shifted" derivative $\tilde \pa_t = \pa_t + \frac{2}{\h}\,{\rm id}$ instead of
$\pa_t$. The partial derivatives remain commutative, while the permutation
relations become
\be
\begin{array}{l@{\quad}l@{\quad}l@{\quad}l}
\tilde\pa_t\,t - t\,\tilde\pa_t = \hh\,\tilde\pa_t & \tilde\pa_t\, x - x\,\tilde\pa_t
=-\hh\,\pa_x &
\tilde\pa_t\, y - y\, \tilde\pa_t=-\hh\,\pa_y &\tilde\pa_t\, z - z\,\tilde\pa_t=- \hh\,\pa_z\\
\rule{0pt}{7mm}
\pa_x\, t - t\,\pa_x = \hh\,\pa_x &\pa_x \,x -  x\,\pa_x = \hh\,\tilde\pa_t &
\pa_x \, y-  y\,\pa_x = \hh\,\pa_z & \pa_x \,z - z\, \pa_x  = - \hh\,\pa_y \\
\rule{0pt}{7mm}
\pa_y \,t - t \, \pa_y = \hh\,\pa_y & \pa_y \,x -  x\,  \pa_y = -\hh\,\pa_z &
\pa_y \,y - y \,  \pa_y = \hh\,\tilde\pa_t & \pa_y \,z - z \,  \pa_ y= \hh\,\pa_x\\
\rule{0pt}{7mm}
\pa_z \,t - t \,\pa_z = \hh\,\pa_z & \pa_z \,x - x \,\pa_z = \hh\,\pa_y&
\pa_z \,y -  y\,\pa_z = -\hh\,\pa_x & \pa_z \,z - z \,\pa_z = \hh\,\tilde\pa_t.
\end{array}
\label{leib-r}
\ee

In the generators $t,\,x,\,y,\,z,\,\tilde\pa_t,\,\pa_x,\, \pa_z,\,\pa_z$,
the algebra $\W(U(u(2)_\h))$ can be treated as the enveloping algebra of a
semi-direct product of the commutative Lie algebra generated by the partial
derivatives $\tilde\pa_t,\,\pa_x,\, \pa_z,\,\pa_z$  and that $u(2)_\h $:  the latter
algebra acts onto the former one in accordance with formulae (\ref{leib-r}).
Then, from the Poincar\'e-Birkhoff-Witt theorem  it follows that the graded
algebra $\Gr\,\W(U(u(2)_\h))$ associated with the Weyl algebra
$\W(U(u(2)_\h))$ is canonically isomorphic to the commutative algebra
$\Sym(W)$ of the vector space $W=\span(t,\,x,\,y,\,z,\pa_t,\,\pa_x,\,
\pa_z,\,\pa_z).$

It is not difficult to compute the Poisson brackets arising on the algebra
$\Sym(W)$
in the limit  $\h\rightarrow 0$ of the algebra $\W(U(gl(2)_\h)$ (with the generator
$\tilde\pa_t$): it is sufficient to set $\h=1$ in  formulae (\ref{leib-r}). We denote
this Poisson bracket $\{\,,\,\}_1$. Thus we have
$$\{{p}_t, t\}_1={{{p}_t}\over 2},\quad \{p_x, t\}_1={{{p}_t}\over 2}
\quad{\rm etc}. $$

Here we replaced the symbols of the partial derivatives
$\{\pa_t,\,\pa_x,\,\pa_y,\,\pa_z\}$ by
the corresponding momenta $\{p_t,\,p_x,\,p_y,\,p_z\}$.

Besides, on the algebra  $\Sym(W)$ there exists the usual Darboux bracket
(denoted
$\{\,,\,\}_0$):
$$
\{{p}_t, t\}_0=1,\quad \{{p}_t, x\}_0=0 \quad{\rm etc}.
$$

The brackets $\{\,,\,\}_i$, $i=0,1$, are compatible with each other. The algebra
$\W(U(gl(2)_\h)$
(with  the generator $\pa_t$) can be treated as a quantization of their sum.

One more bracket on the algebra $\Sym(W)$ is the semi-classical counterpart  of
the algebra $\W(\M)$. We denote this bracket $\{\,,\,\}_2$.
In the coordinates $\{l=a+d,\,h=a-d,\,b,\,c\}$ and the corresponding momenta
$\{p_l,\,p_h,\,p_b,\,p_c\}$ the Poisson structure $\{\,,\,\}_2$ is given by the
following table
$$
\begin{array}{lll}
\{h,b\}_2=-2b(h+l)\quad &\{h,c\}_2=2c(h+l)\quad&\{b,c\}_2=-h(h+l)\\
\rule{0pt}{4.5mm}
 \{l,h\}_2=0 & \{l,b\}_2=0& \{l,b\}_2= 0\\
 \rule{0pt}{6mm}
 \{p_h,p_b\}_2 = 2p_b(p_h-p_l)&\{p_h,p_c\}_2 = -2p_c(p_h-p_l)& \{p_b,p_c\}_2
 =4p_h(p_h-p_l)\\
 \rule{0pt}{4.5mm}
 \{p_l,p_h\}_2 = 0 &  \{p_l,p_b\}_2 = 0 & \{p_l,p_c\}_2 =0
 \end{array}
$$
$$
\begin{array}{ll}
\rule{0pt}{7mm}
\{l,p_l\}_2 = 2+lp_l+hp_h+bp_b+cp_c\qquad &\{h,p_l\}_2 = lp_h+hp_l+bp_b-
cp_c   \\
\rule{0pt}{5mm}
\displaystyle \{b,p_l\}_2 = b(p_h-p_l)+\frac{1}{2}\,p_c(h+l)&\displaystyle
\{c,p_l\}_2 = c(p_h+p_l)+ \frac{1}{2}\,p_b(l-h)\\
\rule{0pt}{6mm}
\{l,p_h\}_2 = lp_h+hp_l-bp_b+cp_c & \{h,p_h\}_2 = 2+hp_h+lp_l+3bp_b-cp_c\\
\rule{0pt}{5mm}
\{b,p_h\}_2 = -b(p_h-p_l)+\displaystyle\frac{1}{2}\,p_c(h+l)&\displaystyle
\{c,p_h\}_2 = c(p_h+p_l)-\frac{1}{2}\,lp_b-\frac{3}{2}\,hp_b\\
\rule{0pt}{6mm}
\{l,p_b\}_2 = p_b(h+l) -2c(p_h-p_l)& \{h,p_b\}_2 = p_b(h+l)+2c(p_h-p_l)   \\
\rule{0pt}{5mm}
\{b,p_b\}_2 =2 + 2bp_b+(h+l)(p_h-p_l) & \{c,p_b\}_2 = 0  \\
\rule{0pt}{6mm}
\{l,p_c\}_2 = p_c(l-h)+2b(p_l+p_h)& \{h,p_c\}_2 = 2b(p_l-3p_h)-p_c(l-h)\\
\rule{0pt}{5mm}
\{b,p_c\}_2 = 0& \{c,p_c\}_2 = 2+2cp_c+l(p_l+p_h)-hp_l+3hp_h.
\end{array}
$$
Note, that on the complexification of the algebra $\Sym(W)$ the bracket $\{\,,\}_2$ is compatible with those $\{\,,\}_0$ and $\{\,,\}
_1$. Thus, we get a pencil generated by three Poisson brackets $\{\,,\}_i,\,\,i=0,1,2$.
Whereas,  the algebra $\W(\N)$ corresponding to the standard Hecke symmetry is a
quantum counterpart of this Poisson pencil.

Now, let us turn to the permutation relations (\ref{leib-r}). They enable us to
compute
similar  relations between the partial derivatives and any power of a generator
of the algebra $U(u(2)_\h)$.

\begin{proposition}
\label{prop:6}
For all $k=0,1,2...$ the following relations hold true
$$
\begin{array}{ll}
\displaystyle
\dd\,t^k=\Big(t+{\h\over 2}\Big)^k\,\dd & \dd\, x^k= A_k(x)\, \dd-B_k(x)\,\pa_x\\
\rule{0pt}{5mm}
\dd\, y^k= A_k(y)\, \dd-B_k(y)\,\pa_y \qquad & \dd\, z^k= A_k(z)\, \dd-B_k(z)
\,\pa_z\\
\rule{0pt}{7mm}
\displaystyle
\pa_x\,t^k=\Big(t+{\h\over 2}\Big)^k\,\pa_x &\pa_x\, x^k=A_k(x)\,\pa_x+ B_k(x)\,
\dd\\
\rule{0pt}{5mm}
\pa_x\, y^k=A_k(y)\,\pa_x+ B_k(y)\, \pa_z & \pa_x\, z^k= A_k(z)\, \pa_x-B_k(z)
\,\pa_y\\
\rule{0pt}{7mm}\displaystyle
\pa_y\,t^k=\Big(t+{\h\over 2}\Big)^k\,\pa_y & \pa_y\, x^k= A_k(x)\, \pa_y-B_k(x)
\,\pa_z\\
\rule{0pt}{5mm}
\pa_y\, y^k=A_k(y)\,\pa_y+ B_k(y)\, \dd & \pa_y\, z^k=A_k(z)\,\pa_y+ B_k(z)\,
\pa_x\\
\rule{0pt}{7mm}\displaystyle
\pa_z\,t^k=(t+{\h\over 2})^k\,\pa_z & \pa_z\, x^k=A_k(x)\,\pa_z+ B_k(x)
\,\pa_y\\
\rule{0pt}{5mm}
\pa_z\, y^k= A_k(y)\, \pa_z-B_k(y)\,\pa_x & \pa_z\, z^k=A_k(z)\,\pa_z+ B_k(z)\,
\dd,
\end{array}
$$
where
$$
A_k(v)=\frac{1}{2}\,\bigg(\Big(v-i\,{\h\over 2}\Big)^k+\Big(v+i\,{\h\over 2}\Big)
^k\bigg),\qquad B_k(v)=
\frac{i}{2}\,\bigg(\Big(v-i\,{\h\over 2}\Big)^k-\Big(v+i\,{\h\over 2}\Big)^k\bigg).
$$
\end{proposition}

{\bf Proof.} The proof can be done by induction in the power $k$ of a
monomial.\hfill\rule{6.5pt}{6.5pt}
\medskip

The above formulae enable us to compute the permutation relations between
partial
derivatives and any polynomial or a formal series in one variable.

\begin{corollary}
Let  $f(v)$ be a polynomial or a formal series in one variable $v$ and $A$ and
$B$ be difference operators defined by
$$
A(f(v))={{f(v-i\,{\h\over 2})+f(v+i\,{\h\over 2})}\over 2}\,,\qquad
B(f(v))={{i\,f(v-i\,{\h\over 2})-i\,f(v+i\,{\h\over
2})}\over 2}\,.
$$
Then all formulae in Proposition \ref{prop:6} remain valid if we respectively
replace the powers $v^k$ where $v\in \{t,\,x,\,y,\,z\}$ for $f(v)$ and the factors
$(v+{\h\over 2})^k$, $A_k(v)$ and $B_k(v)$  by $f(v+{\h\over 2})$, $A(f(v))$
and $B(f(v))$.
\end{corollary}

The permutation relations described above enable us to find the result $\pa_v(f)$
of applying a partial derivative $\pa_v$, $v\in \{t,x,y,z\}$, to any polynomial or a
formal
series
$f$. According to our general recipe (\ref{action-ch}), we should move the symbol
of
the partial derivative to the most right position and then apply the relation
$$
\pa_v(1_{\cal W})=0
$$
motivated by (\ref{cu}). For example, we can find
$$
\pa_x(x^k) = \frac{2}{\h}\,B_k(x),
$$
since $\tilde\pa_t(1_{\cal W}) = 2/\h$, etc.

Now, using this method we compute the result of applying a partial derivative
to a decomposable element $f=f_0(t)\,f_1(x)\,f_2(y)\,f_3(z)$ of the algebra
$U(u(2)_\h)$
or its completion $U(u(2)_\h)[[\h]]$. Here $f_i$ are polynomials or formal
series in
one variable.

\begin{proposition}
The following relations hold true for the element $f=f_0(t)\,f_1(x)\,f_2(y)\,f_3(z)$
$$
\dd(f)=2\hhh f_0(t+{\h\over 2})\,(A(f_1)\,A(f_2)\,A(f_3)-B(f_1)\,B(f_2)\,B(f_3)),
$$
$$\pa_x(f)=2\hhh f_0(t+{\h\over 2})\,(B(f_1)\,A(f_2)\,A(f_3)+A(f_1)\,B(f_2)
\,B(f_3)), $$
$$\pa_y(f)=2\hhh f_0(t+{\h\over 2})\,(A(f_1)\,B(f_2)\,A(f_3)-B(f_1)\,A(f_2)\,B(f_3)),
$$
$$\pa_z(f)=2\hhh f_0(t+{\h\over 2})\,(A(f_1)\,A(f_2)\,B(f_3)+B(f_1)\,B(f_2)
\,A(f_3)).
$$
Consequently,
$$
\pa_t(f)=2\hhh \Big(f_0(t+{\h\over 2})\,(A(f_1)\,A(f_2)\,A(f_3)-B(f_1)\,B(f_2)\,
B(f_3))- f\Big).
$$
\end{proposition}

It is interesting to find "eigenfunctions" of all partial derivatives. In the classical
case such eigenfunctions are exponential functions. Consider an element
\be
f_{\Al}=exp(\al_0\, t)\,exp(\al_1\, x)\,exp(\al_2\, y)\,exp(\al_3\, z),\qquad \Al=
(\al_0,\,\al_1,\,\al_2,\,\al_3)\in \h\,\K^4,
\label{exp}
\ee
belonging to the completion $U(u(2)_\h)[[\h]]$.

By using the last proposition we get the following.
\begin{proposition}
$$
\dd(f_{\Al}) =2\hhh exp(\al_0\,{\h\over 2})\,(\cos(\al_1\,{\h\over 2})
\,\cos(\al_2\,{\h\over 2})\, \cos(\al_3\,{\h\over 2})-
\sin(\al_1\,{\h\over 2})\,\sin(\al_2\,{\h\over 2})\, \sin (\al_3\,{\h\over 2}) )\, f_{\Al},
$$
$$
\pa_x(f_{\Al})=2\hhh exp(\al_0\,{\h\over 2})\,(\sin(\al_1\,{\h\over 2})
\,\cos(\al_2\,{\h\over 2})\, \cos(\al_3\,{\h\over 2})+
\cos(\al_1\,{\h\over 2})\,\sin(\al_2\,{\h\over 2})\, \sin (\al_3\,{\h\over 2}) )\, f_{\Al},
$$
$$
\pa_y(f_{\Al}) =2\hhh exp(\al_0\,{\h\over 2})\,(\cos(\al_1\,{\h\over 2})
\,\sin(\al_2\,{\h\over 2})\, \cos(\al_3\,{\h\over 2})-
\sin(\al_1\,{\h\over 2})\,\cos(\al_2\,{\h\over 2})\, \sin (\al_3\,{\h\over 2}) )\, f_{\Al},
$$
$$\pa_z(f_{\Al}) =2\hhh exp(\al_0\,{\h\over 2})\,(\cos(\al_1\,{\h\over 2})
 \,\cos(\al_2\,{\h\over 2})\, \sin(\al_3\,{\h\over 2})+
\sin(\al_1\,{\h\over 2})\,\sin(\al_2\,{\h\over 2})\, \cos (\al_3\,{\h\over 2}) )\, f_{\Al}.
$$
\end{proposition}

Consequently, in our NC setting the ordered exponents (\ref{exp}) are still
"eigenfunctions" of all partial derivatives but the corresponding eigenvalues
are modified with respect to the classical case.

Now, we are able to define an analog  of the de Rham operator
on the algebra $U(u(2)_\h)$. Let $ \bigwedge(u(2))$ be the usual
skew-symmetric algebra with four generators  $dt,\, dx,\, dy,\, dz$ and
$ \bigwedge^k(u(2))$ be its degree $k$ homogeneous component. Introduce the
space of $k$-differential forms on $U(u(2)_\h)$ by
$ \bigwedge^k(u(2))\ot  U(u(2)_\h)$.
Then define an analog of the de Rham operator on the algebra  $U(u(2)_\h)$   as
follows
$$
d(f)=dt\,\pa_t(f)+ dx\,\pa_x(f)+dy\,\pa_y(f)+dz\,\pa_z(f),
$$
where $f$ be an arbitrary element of this algebra. On the higher order differential
forms the action of $d$ is extended  in the usual way
$$
d(\omega \, f)=\omega\, d(f),\,\,\, \omega\in \bigwedge(u(2)),\, f\in  U(u(2)_\h).
$$

Note, that we do not transpose elements of $U(u(2)_\h)$ and those of $
\bigwedge(u(2))$. Nevertheless, the de Rham operator $d$ is well defined
and due to the commutativity of the partial derivatives it is easy to see that
  $d^2=0$, i.e. $d$ is a differential indeed. Emphasize that this way of defining
  the de Rham operator $d$ does not imply any Leibniz rule for it looking like
  the classical one: $d(f\, g)=d(f) \,g+f\, d(g)$.

In the same way it is possible to define
analogs of differential forms and the de Rham operators for any enveloping
algebra  $U(u(m|n)_\h)$.

\begin{remark}\rm
As for the de Rham complex over the  RE algebras there are different approaches
to its
construction. One of them based on Koszul type complexes was considered in
\cite{GS2}
(also, see the references therein).
Another one was developed for the RE algebra playing the role of
q-Minkowski space in \cite{OSWZ,K,AKR}. This approach involves
some permutation rules between "functions" and "differentials" and the classical form of the Leibniz rule for the de Rham differential.

This approach fails provided the RE algebra $\M$ is replaced by its modified counterpart $\N$.
Hopefully, the method above of constructing the de Rham differential on the algebra $U(u(m|n)_\h)$ can be applied
mutatis mutandis to the algebra $\N$.
We plan to go back to this construction in our subsequent publications.
\end{remark}

\section{Laplace operator on $U(u(2))$ and its radial part}

Consider the algebra $\W(U(u(2)_\h))$ in more detail. This algebra is covariant
with respect to spacial rotations, i.e. actions of elements of the group $SO(3)$
(we treat $x,\,y,\,z$ as
"space generators" and $t$ plays the role of the time).  Indeed,
rewriting the generating matrix (\ref{MMM}) through the generators $t,\, x,\, y,\, z$
\be
N=\left(\begin{array}{cc}
t-i\,z& -i\, x-y\\
-i\, x+y& t+i\, z
\end{array}\right)
\label{matrN}
\ee
we can realize an action of the group $SO(3)$ as a similarity transformation
\be
N\longrightarrow  g^{-1}N\, g,
\label{act}
\ee
where $g$ is an element of the group $SU(2)$ arising from the spinor
representation of the group $SO(3)$ (though this representation is
two-fold it does not affect the
result). In the same way the group $SO(3)$ acts on the matrix of partial
derivatives
(\ref{DDD}) which in terms of $\pa_t,\,\pa_x,\,\pa_y,\,\pa_z$ is  as follows
$$
D= {1\over 2} \left(\begin{array}{cc}
\pa_t+i\,\pa_z& i\, \pa_x+\pa_y\\
i\,\pa_x-\pa_y& \pa_t-i\pa_z
\end{array}\right).
$$

\begin{remark}
\rm
Note that the algebra $\W(U(u(2)_\h))$ is not covariant with respect
to the
boosts\footnote{And we do not know any deformation of the Poincar\'e group in
the spirit of the $\kappa$-Poincar\'e one
(see \cite{MR}).}). However, since the subalgebra $\D\subset\W(U(u(2)_\h))$
formed by partial derivatives is commutative,  it can be
equipped with any metric, for example, Euclidian or Lorentzian ones.
Consequently, we can define
the Laplace operator  of any signature, for instance  $(4,0)$ or
 $(1,3)$ (in the latter case it is called d'Alembertian and is denoted $\square $).
Thus,  we can define  the Klein-Gordon equation
$$
(\square-m^2)\, f= 0,\qquad
\square=\pa_t^2-\pa_x^2-\pa_y^2-\pa_z^2
$$
in the classical way but with $f$ being an element of the algebra $U(u(2)_\h)$ or
its completion and  the  partial derivatives defined above.

In a similar way we can define analogs of the Dirac and Maxwell operators. Namely, we introduce
them by the classical formulae
but with a new meaning of the partial derivatives.
\end{remark}

In what follows we concentrate ourselves to the Laplace operator
$$\De=\pa_x^2+\pa_y^2+\pa_z^2$$
which is invariant with respect to the spacial rotations.
In the remaining part of the paper we define and compute its "radial
part".

Consider
the center $Z=Z(U(u(2)_h))$ of the algebra $U(u(2)_h)$. It is generated by the
elements $\Tr\, N$, $\Tr\, N^2$, or, equivalently, by
$$
t = \frac{\Tr N}{2}\qquad {\rm and}\qquad
\Cas=x^2+y^2+z^2 = \frac{1}{4}\,((\Tr N)^2 - 2 \Tr N^2).
$$
The latter element is the quadratic Casimir of the algebra $U(su(2)_h)$.
More precisely, any element $f\in Z$ is a polynomial in $t$ and $\Cas $. Below
we show that
$$
\forall\,f\in Z\Rightarrow\Delta(f)\in Z.
$$
This fact motivates the following definition.

\begin{definition}\label{def:11}
The operator $\De $ restricted to the center $Z$  is called the radial part
of the Laplacian $\De $ and is denoted $\De_{rad}$.
\end{definition}

In a similar way we can define the radial part of any linear
combination of operators
$t^k\,\Cas^p\,\pa_t^m\,\De^n$ invariant with respect to rotations.

Note that in  the classical case  the radial part of the Laplace operator
on the vector space $\Mat(m)$ is realized via eigenvalues of symmetric
 or Hermitian  matrices.
Fortunately, in the RE algebras (modified or not) and their
$q=1$ counterparts  the generating matrices satisfy an analog of  the
Cayley-Hamilton identity enabling us to introduce  quantum version of
eigenvalues.  Let us exhibit their construction.

The matrix $N$ defined in (\ref{matrN}) satisfies the Cayley-Hamilton (CH) identity
$$
%N^2-(a+d+\h)\, N+(ad-bc+\h\,a)\,I=
N^2-(2\,t+\h)\, N+\,(t^2+x^2+y^2+z^2+\h\, t)\,I= 0,
$$
which can be checked by direct calculations.

Let $\mu_1,\, \mu_2$ be the roots of the equation
$$
\mu^2 - (2\,t+\h)\,\mu + (t^2+x^2+y^2+z^2+\h\, t) = 0.
$$
This means that $\mu_1$ and $\mu_2$ satisfy the relations
$$
\mu_1+\mu_2=2\,t+\h,\qquad
\mu_1 \,\mu_2= t^2+x^2+y^2+z^2+\h\, t.
$$
Since the coefficients of the CH identity are central, $\mu_1$ and $\mu_2$
belong to the algebraic extension of the center $Z$.

Then the quantity $\Tr N^k$, $k\in {\Bbb Z}_+$, can be expressed via these roots
by
the formula\footnote{This formula can be obtained from its $q$-analog
established in
\cite{GS1} in the limit $q \to 1$.}
$$
\Tr N^k= {{\mu_1-\mu_2-\h}\over{\mu_1-\mu_2}}\, \mu_1^k+{{\mu_2-\mu_1-\h}
\over{\mu_2-\mu_1}}\, \mu_2^k.
$$
In particular,
$$
\Tr N = \mu_1+\mu_2-\h,\qquad
\Tr N^2= \mu_1^2+\mu_2^2-\h\,(\mu_1+\mu_2).
$$
Since, on the other hand, $\Tr N = 2t$ and $\Tr N^2 = 2(t^2- \Cas)$, we get
\be
t={{\mu_1+\mu_2-\h}\over 2},\qquad
\Cas = -{{(\mu_1-\mu_2)^2-\h^2}\over 4}.
\label{mu-param}
\ee
So, any element of $Z$, being a polynomial in $t$ and $\Cas $, can be
expressed as a {\it symmetric
polynomial} in the roots $\mu_1$ and $\mu_2$.

Our next step consists in computing the quantities
$$
\De(t^k\, \Cas^p),\quad k,\, p\in {\Bbb Z}_+
$$
and presenting  the result through  the roots $\mu_1,\, \mu_2$.
To this end we introduce a first order operator $Q=x\, \pa_x+y\,\pa_y+z\,\pa_z$
and consider four second order operators
$$
\De_0=\dd^2,\quad \De_1=\De=\pa_x^2+\pa_y^2+\pa_z^2,\quad \De_2=Q\,
\dd,\quad \De_3=
Q^2.
$$
Our immediate aim is to compute permutation relations of all these operators and
the elements $t^k$ and $\Cas^p$.

It is not difficult to see that
\be
\De_i\, t^k= (t+\h)^k\, \De_i,\quad \forall\, i=0,1,2,3.
\label{act-t}
\ee
Thus, we have only to compute the permutation relations of the operators
$\Delta_i$
and the Casimir element $\Cas$.

\begin{proposition}
The following permutation relations hold true:
$$
\Delta_i \Cas = \sum_{j=0}^3\Pi_{ij}\,\Delta_j, \quad 0\le i\le 3,
$$
where the matrix $\Pi = \|\Pi_{ij}\|$ reads
$$
\Pi =
\left(\begin{array}{c@{\hspace*{6.5mm}}c@{\hspace*{6.5mm}}
c@{\hspace*{6.5mm}}c}
\displaystyle\Cas-{3\over 2}\,\h^2 &\displaystyle{\h^2\over 2} &-2\h &0 \\
\rule{0pt}{8mm}
\displaystyle {3\over 2}\,\h^2 &\displaystyle \Cas-{\h^2\over 2} &2\h& 0\\
\rule{0pt}{8mm}
\h\,\Cas &0 &\displaystyle\Cas-{\h^2\over 2} &-\h \\
\rule{0pt}{8mm}
\h^2\Cas&\displaystyle -{\h^2\over 2}\,\Cas &
\displaystyle\h\left(2\Cas+{\h^2\over 4}\right)&
\displaystyle\Cas+{\h^2\over 2}
\end{array}\right).
$$
\end{proposition}

{\bf Proof} can be done by straightforward computations on the base of
(\ref{leib-r}).
\hfill\rule{6.5pt}{6.5pt}
\medskip

Now, we immediately get
$$
\Delta_i\Cas^p = \sum_{j=0}^3(\Pi^p)_{ij}\Delta_j.
$$
In order to get the action of the operators $\Delta_i$ we apply the rule of action
of the partial derivatives on the unit element $1_{\W}$
$$
\De_0(1_{\W})=4\h^{-2},\quad \De_i(1_{\W})=0,\quad i=1,2,3.
$$
So, we arrive to  the following result
\be
\Delta_i(\Cas^p) = 4\h^{-2}(\Pi^p)_{i0}.
\label{del-cas}
\ee

To calculate the $p$-th power of the matrix $\Pi $ we substitute the
parametrization (\ref{mu-param})
for ${\rm Cas}$ and calculate the spectrum of $\Pi $. By a direct calculation one
can verify that the matrix $\Pi $
is  semisimple:
$$
\Pi\sim {\rm diag}(\lambda_0,\lambda_0,\lambda_+,\lambda_-),
$$
where the eigenvalues are
$$
\lambda_0= \frac{1}{4}\,(\h^2 - (\mu_1-\mu_2)^2),\quad
\lambda_\pm = \frac{1}{4}\, (\h^2-(\mu_1-\mu_2\pm 2\,\h)^2)\,.
$$

Now, the matrix value $f(\Pi)$ of a function $f(x)$ (provided that $f(x)$ is defined
on
the spectrum
of $\Pi $) can be found in terms of the Lagrange-Sylvester polynomial (see
\cite{G})
\be
f(\Pi) = \sum_{a=0,\pm}f(\lambda_a)\prod_{b=0,\pm\atop b\not=a}
\frac{(\Pi-\lambda_bI)}{(\lambda_a-\lambda_b)}.
\label{LS-pol}
\ee
Taking, in particular, $f(x) = x^p$ we get the following explicit formulae from
(\ref{del-cas})
$$
\begin{array}{l}
\displaystyle\Delta_0(\Cas^p) = \frac{1}{\h^2(\mu_1-\mu_2)}\,\Bigl(
2(\mu_1-\mu_2)\lambda_0^p+(\mu_1-\mu_2+2\h)\lambda_+^p +
(\mu_1-\mu_2-2\h)\lambda_-^p\Bigr),\\
\rule{0pt}{7mm}
\displaystyle\Delta_1(\Cas^p) = \frac{1}{\h^2(\mu_1-\mu_2)}\,\Bigl(
2(\mu_1-\mu_2)\lambda_0^p - (\mu_1-\mu_2+2\h)\lambda_+^p -
(\mu_1-\mu_2 -2\h)\lambda_-^p\Bigr),\\
\rule{0pt}{7mm}
\displaystyle\Delta_2(\Cas^p) = \frac{1}{2\h^2(\mu_1-\mu_2)}\,\Bigl(
-2\h(\mu_1-\mu_2)\lambda_0^p +((\mu_1-\mu_2)^2+\h(\mu_1-\mu_2)-2\h^2)
\lambda_+^p\\
\hspace*{50mm} -((\mu_1-\mu_2)^2 -\h(\mu_1-\mu_2)-2\h^2)\lambda_-^p\Bigr),\\
\displaystyle\Delta_3(\Cas^p) = \frac{1}{4\h^2}\,\Bigl(
(2\h^2-(\mu_1-\mu_2)^2)\lambda3_0^p +((\mu_1-\mu_2)
^2+\h(\mu_1-\mu_2)-2\h^2)\lambda_+^p\\
\hspace*{50mm} + ((\mu_1-\mu_2)^2 -\h(\mu_1-\mu_2)-2\h^2)\lambda_-^p\Bigr).
\end{array}
$$

The second formula of this list allows us to explicitly calculate  the radial part
$\Delta_{rad}$
(recall that $\De=\De_1$), the radial parts of  other operators above can be
calculated in the same way.
Indeed, as follows from Definition \ref{def:11}, we have to know the result of
action of the operator $\Delta_1$ on any symmetric function in $\mu_1$ and
$\mu_2$. It is convenient to fix the following set of generators for the ring of
such functions:
$$
\mu_1+\mu_2\qquad{\rm and}\qquad (\mu_1-\mu_2)^2 = \h^2 - 4\Cas.
$$

Now, extend the second formula from the list above to any function $g(\Cas)$
where $g$ is a polynomial or a formal series in one variable. We have
$$
\begin{array}{l}
\displaystyle\Delta(g(\Cas)) = \frac{1}{\h^2(\mu_1-\mu_2)}\,\Bigl(
2(\mu_1-\mu_2)g(\lambda_0) - (\mu_1-\mu_2+2\h)g(\lambda_+) -
(\mu_1-\mu_2 -2\h)g(\lambda_-)\Bigr).
\end{array}$$

Then, by applying this formula to the function  $f(\h^2-4\Cas)$ and using (\ref{del-cas})
and (\ref{LS-pol}) we finally get
\begin{eqnarray}
\Delta_{rad}\left(f((\mu_1-\mu_2)^2)\right)& =& \frac{1}{\h^2}
\,\Bigl(2f((\mu_1-\mu_2)
^2) -f((\mu_1-\mu_2-2\h)^2)-f((\mu_1-\mu_2+2\h)^2)\Bigr)\nonumber\\
&& +\frac{2}{\mu_1-\mu_2}\,\frac{1}{\h}\,\left(f((\mu_1-\mu_2-2\h)^2)
-f((\mu_1-\mu_2+2\h)^2)\right).\label{rad-p}
\end{eqnarray}

The action on a function depending on the both generators $\mu_1+\mu_2$
and $(\mu_1-\mu_2)^2$ are easily obtained with the help of (\ref{act-t}).
In order to better visualize the final formula we set $\la=\mu_1+\mu_2$ and
$\mu=(\mu_1-\mu_2)^2$.
Then we have
\begin{eqnarray*}
\Delta_{rad}\left(f(\la,\mu)\right)& =& \frac{1}{\h^2}\,\Bigl(2f(\la+2\h, \mu)
 -f(\la+2\h, \mu+4\h^2+4\h \sqrt{\mu})-f(\la+2\h, \mu+4\h^2-4\h \sqrt{\mu})\Bigr)
 \nonumber\\
&& +\frac{2}{\sqrt{\mu}}\,\frac{1}{\h}\,\left(f(\la+2\h, \mu+4\h^2-4\h \sqrt{\mu})
-f(\la+2\h, \mu+4\h^2+4\h \sqrt{\mu})\right).
\end{eqnarray*}
Note that although the function $\mu\to \sqrt{\mu}$ is twofold it does not affect the
result
(here we are dealing with functions $f(\lambda,\mu)$ which are polynomials or
formal series
in commutative variables $\lambda$ and $\mu$).

In the limit $\h\to 0$ the difference operator $\Delta_{rad}$ turns into the following
second
order differential operator
$$
-16\mu {{\pa^2}\over {\pa\mu^2}}-24 {{\pa}\over {\pa\mu}}.
$$
Being rewritten via the variable $r$ such that $\mu=-4 r^2$ we get the usual radial
part of the
classical Laplacian on ${\R}^3$
$$
{{\pa^2}\over {\pa r^2}}+{2\over r}{{\pa}\over {\pa r}}.
$$

By completing the paper we want to design a  plan  of a possible application of
our
method to constructing two-parameter deformations of  the rational Calogero-Moser models.

Let $\N$ be the standard mRE algebra (it is a braided deformation of the enveloping
algebra $U(gl(m)_\h)$). Also, let
$D$ be the matrix of the partial derivatives on this mRE algebra as introduced in the section 2.
Consider the  operators $ Tr_{R^{-1}} D^k, k=0,1,2,..., m$ acting on the algebra $\N$.
(Recall that $D$ is the generating matrix of the RE algebra corresponding to the Hecke symmetry ${R^{-1}}$.)
They  commute with each other. Besides, they map the center
$Z=Z(U(gl(m)_\h))$ of the algebra $U(gl(m)_\h)$ into itself.
 Consequently, the restriction of the operators $ Tr_{R^{-1}} D^k$ to $Z$ is well
 defined. By expressing these restricted operators
 via the eigenvalues of the generating matrix $N$ of the algebra $\N$ we get a family
 of operators in involution. Hopefully, these operators are difference ones
 and they are two-parameter deformations of the corresponding classical
 differential operators which are gauge equivalent to the rational
 Calogero-Moser operator and its higher counterparts respectively. However, computations in
 higher dimensional case become much harder.

 It would be also interesting to apply this method to   algebras which are
 deformations of the super-algebra  $U(gl(m|n)_\h)$. Note that radial parts of
 Laplace operators on supergroups have been studied by F.Berezin (see
 \cite{B}).


\begin{thebibliography}{IOPPP}

\bibitem[AKR]{AKR} de Azc\'arraga J., Kulish P., R\'odenas F. {\em Reflection equation and
$q$-Minkowski space algebras}, Lett. Math. Phys. 32 (1994) no.3, 173--182.

\bibitem[B]{B} Berezin F. {\em Introduction to super-analysis}, D.Reider
publishing company 1987.

\bibitem[FP]{FP} Faddeev L., Pyatov P. {\em The Differential calculus on
quantum linear groups}, Amer. Math. Soc. Transl. Ser. 2, 175 (1996) 35--47.

\bibitem[FRT]{FRT}
Faddeev L., Reshetikhin N., Takhtadzhyan L. {\em Quantization of Lie
Groups and Lie Algebras}, English translation: Leningrad Math J. 1 (1990)
193--225.

\bibitem[G]{G} Gantmacher {\em The theory of Matrices}, vol. 1, AMS Chelsea
publ., Providence, Rode Island, 2000.

\bibitem[GR]{GR} Gelfand I., Retakh V. {\em Quansideterminants. I.}, Selecta
Math. (N.S.) 3 (1997), 517--546.

\bibitem[GPS1]{GPS1}
Gurevich D., Pyatov P., Saponov P. {\em The Cayley-Hamilton
theorem
for quantum matrix algebras of $GL(m|n)$ type}, English
translation:
St.Petersburg Math.J. 17 (2006) 119--135.

\bibitem[GPS2]{GPS2}
Gurevich D., Pyatov P., Saponov P. {\it Representation theory of
(modified) Reflection Equation Algebra of the $GL(m|n)$ type}, St.Petersburg
Math.J.
English translation:  20 (2009) 213--253.

\bibitem[GPS3]{GPS3} Gurevich D., Pyatov P., Saponov P. {\it Braided
differential operators on quantum algebras},
J. of Geometry and Physics 61 (2011) 1485--1501.

\bibitem[GS1]{GS1} Gurevich D., Saponov P. {\it Quantum line bundles via
Cayley-Hamilton identity} J. Phys. A 34 (2001)  4553--4569.

\bibitem[GS2]{GS2}  Gurevich D., Saponov P.  {\em Braided affine geometry and
$q$-analogs of wave operators},
J. Phys. A 42 (2009) no. 31, 51 pp.

\bibitem[IOP]{IOP} Isaev A., Ogievetsky O., Pyatov P.
{\em On quantum matrix algebras satisfying the
Cayley-Hamilton-Newton identities}, J.Phys. A: Math. Gen. 32
(1999) L115--L121.

\bibitem[IP]{IP} Isaev A., Pyatov P. {\em Covariant differential calculus
complexes on quantum linear groups}, J.Phys. A 28
(1995) 8, 2227--2246.

\bibitem[IV]{IV} Isaev A., Vladimirov A. {\em $GL_q(N)$-covariant Braided Differential Bialgebras},
Lett. Math. Physics 33 (1995), 297--302.

\bibitem[K]{K} Kulish P.P. {\em Representations of $q$-Minkowski space algebra},
Algebra i Analiz 6 (1994), no. 2, 195--205 (in Russian).\par

English transl.: St. Petersburg Math. J. 6 (1995), no.2, 365--374.

\bibitem[M]{M} Majid S. {\em Foundations of quantum group
theory}, Cambridge University Press, Cambridge, 1995.

\bibitem[MR]{MR} Majid S., Ruegg H., {\em Bicrossproduct structure of $\kappa$-
Poincare' group and non-commutative geometry},
Phys. Lett. B 334 (1994) 348--354.

\bibitem[OSWZ]{OSWZ} Ogievetsky O.,  Schmidke W.B., Wess J.,  Zumino B.,
{\em  $q$-deformed Poincare' algebra}, Commun. Math.
Phys. 150, (1992) 495--518.

\bibitem[WZ]{WZ} Wess J., Zumino B. {\em Covariant differential calculus on the
quantum hyperplane}
Nuclear Phys. B Proc. Suppl. 18B (1990), 302--312 (1991)

\bibitem[W]{W} Woronowicz S. {\em Differential Calculus on
Compact Matrix Pseudogroups (Quantum groups)}, CMP 122  (1989),
125--170.

\end{thebibliography}
\end{document}